
\magnification 1200
\parskip 9pt plus 5pt minus 3pt
\def\item{\parskip 5pt plus 3pt minus 3pt \par\hang\textindent}

\headline{\ifodd\pageno\ifnum\pageno=1\else\hfil\rlap{\hautdepage\folio}\fi\else\llap{\hautdepage\folio}\hfil\fi} 
\footline{}              
\catcode`\@=11
\newdimen\margereliured \margereliured=1truecm
\newdimen\margereliureg \margereliureg=-0.1truecm
\def\plainoutput{%
\shipout\vbox{
\ifodd\pageno \moveright\margereliured
              \else \moveleft\margereliureg \fi
\hbox{\vbox{\makeheadline\pagebody\makefootline}}}%
\advancepageno
\ifnum\outputpenalty>-\@MM \else\dosupereject\fi}
\catcode`\@=12

\def\tenpoint{%
  \initfileref
  \textfont0=\tenrm \scriptfont0=\sevenrm \scriptscriptfont0=\fiverm
  \def\rm{\fam0\tenrm}%
  \textfont1=\teni \scriptfont1=\seveni \scriptscriptfont1=\fivei
  \def\oldstyle{\fam1\teni}%
  \textfont2=\tensy \scriptfont2=\sevensy \scriptscriptfont2=\fivesy
  \textfont\itfam=\tenit
  \def\it{\fam\itfam\tenit}%
  \def\sl{\fam\slfam\tensl}%
  \textfont\slfam=\tensl \scriptfont\slfam=\sevensl \scriptscriptfont\slfam=\fivesl
  \def\bf{\fam\bffam\tenbf}%
  \textfont\bffam=\tenbf \scriptfont\bffam=\sevenbf
  \scriptscriptfont\bffam=\fivebf  
  \def\tt{\fam\ttfam\tentt}%
  \textfont\ttfam=\tentt
  \abovedisplayskip=6pt plus 2pt minus 6pt
  \abovedisplayshortskip=0pt plus 3pt
  \belowdisplayskip=6pt plus 2pt minus 6pt
  \belowdisplayshortskip=7pt plus 3pt minus 4pt
  \smallskipamount=3pt plus 1pt minus 1pt
  \medskipamount=6pt plus 2pt minus 2pt
  \bigskipamount=12pt plus 4pt minus 4pt
  \normalbaselineskip=12pt
  \setbox\strutbox=\hbox{\vrule height8.5pt depth3.5pt width0pt}%
  \normalbaselines\rm}
\catcode`\|=13
\def\today{\ifcase\month\or january \or february \or march \or april
\or may \or june\or july\or august \or september\or october\or november\or
December\fi\ \number\day , \number\year}


\newskip\LastSkip
\def\nobreakatskip{\relax\ifhmode\ifdim\lastskip>0pt
  \LastSkip\lastskip\unskip
  \nobreak\hskip\LastSkip
  \fi\fi}
\catcode`\;=\active \def;{\nobreakatskip\string;}
\catcode`\:=\active \def:{\nobreakatskip\string:}
\catcode`\!=\active \def!{\nobreakatskip\string!}
\catcode`\?=\active \def?{\nobreakatskip\string?}

\newif\ifrefvis
\refvisfalse

\newif\ifarxiv
\arxivfalse
\def\arXiv{\arxivtrue}

\newif\iffrance

\def\anglais{\francefalse}

\newskip\afterskip
\catcode`\@=11
\def\p@int{.\par\vskip\afterskip\penalty100} 
\def\p@intir{\discretionary{.}{}{.\kern.35em---\kern.7em}}
\def\pointir{\afterassignment\pointir@\global\let\next=}
\def\pointir@{\ifx\next\par\p@int\else\p@intir\fi\egroup\next}
\catcode`\@=12
\def|{\relax\ifmmode\vert\else\findef\fi}
\def\findef{\errhelp{Cette barre verticale ne correspond ni a un \vert mathematique
                        ni a une fin de definition, le contexte doit vous indiquer ce qui manque.
                        Si vous vouliez inserer un long tiret, le codage recommande est ---,
                        dans tous les cas, la barre fautive a ete supprimee.}%
                        \errmessage{Une barre verticale a ete trouvee en mode texte}}

\def\TITR#1|{\null{\mss\baselineskip=17pt
                           \vskip 3.25ex plus 1ex minus .2ex
                           \leftskip=0pt plus \hsize
                           \rightskip=\leftskip
                           \parfillskip=0pt
                           \noindent #1
                           \par\vskip 2.3ex plus .2ex}}
 
\def\auteur#1|{\penalty 500
               \vbox{\centerline{
                 \iffrance par \else by \fi #1}
                \vskip 10pt}\penalty 500}


\def\fonction#1|#2|#3|#4|{\lower 6pt\hbox{
$\matrix{#1 \hfill &\longrightarrow \hfill & #2\hfill\cr
\hfill #3  &\longmapsto \hfill &#4 \hfill\cr}$}}
\newcount\thesection
\newcount\thesubsection
\newcount\thesubsubsection
\newcount\theparagraf
\newcount\thetheo
\newcount\theequ
\global\thesection=0
\global\thesubsection=0
\global\thesubsubsection=0
\global\theparagraf=0
\global\thetheo=0
\global\theequ=0

\font\sevensl= cmsl7
\font\fivesl= cmsl5

\font\tentite = cmbx10 at 16pt
\font\seventite = cmbx7 at 11pt
\font\fivetite = cmbx5 at 8pt
\newfam\titefam
\textfont\titefam = \tentite
\scriptfont\titefam = \seventite
\scriptscriptfont\titefam = \fivetite

\font\sectfont = cmbx10 at 14pt
\font\sectscript = cmbx7 at 10pt
\font\sectsscript = cmbx5 at 7pt
\newfam\sectionfam 
\textfont\sectionfam = \sectfont 
\scriptfont\sectionfam = \sectscript
\scriptscriptfont\sectionfam = \sectsscript
\def\sectionfont{\fam\sectionfam\sectfont}

\font\subsectfont =  cmbx10 at 12pt
\font\subsectscript = cmbx7 at 8pt
\font\subsectsscript = cmbx5 at 6pt
\newfam\subsectionfam
\textfont\subsectionfam = \subsectfont
\scriptfont\subsectionfam = \subsectscript
\scriptscriptfont\subsectionfam = \subsectsscript
\def\subsectionfont{\fam\subsectionfam\subsectfont}

\font\mss=cmss12 scaled \magstep1

\font\htdpfont =  cmss10
\font\subhtdpfont = cmss9  
\font\ssubhtdpfont = cmss8 
\newfam\htdepagefam
\textfont\htdepagefam = \htdpfont 
\scriptfont\htdepagefam = \subhtdpfont
\scriptscriptfont\htdepagefam = \ssubhtdpfont
\def\hautdepage{\fam\subsectionfam\htdpfont}
\font\sf =cmss10
\font\ding = pzdr at 10pt

\font\tenmsb=msbm10
\font\sevenmsb=msbm7
\font\fivemsb=msbm5
\newfam\msbfam
\textfont\msbfam=\tenmsb
\scriptfont\msbfam=\sevenmsb
\scriptscriptfont\msbfam=\fivemsb
\def\Bbb#1{{\fam\msbfam\relax#1}}

\font\teneufrak=eufm10
\font\seveneufrak=eufm7
\font\fiveeufrak=eufm5
\newfam\eufrak
\textfont\eufrak=\teneufrak
\scriptfont\eufrak=\seveneufrak
\scriptscriptfont\eufrak=\fiveeufrak
\def\mathfrak#1{{\fam\eufrak\relax#1}}


\def\begincentered{\par\begingroup
\def \par{\hss\egroup\line\bgroup\hss}\obeylines
\line\bgroup\hss}
\def\endcentered{\hss\egroup\endgroup}

\hsize=12.5cm
\vsize=20cm
\parindent=1cm
\baselineskip=13pt
\hoffset=-0.1cm
\voffset=0.5cm

\long\def\partie#1{\begingroup\sectionfont
        \par\penalty -500
        \vskip 3.25ex plus 1ex minus .2ex
        \skip\afterskip=1.5ex plus .2ex
        \baselineskip=17pt        
        \par
        \def \par{\hss\egroup\line\bgroup\hss}\obeylines
        \line\bgroup\hss
\global\advance\thesection by 1 
\xdef \lastref{\number\thesection}
\global\thesubsection=0 \global\theparagraf=0  \number\thesection\quad#1\hss\egroup\endgroup\par}

\long\def\subpartie#1{\begingroup\subsectionfont%
                          \par\penalty -200
                          \vskip 3.25ex plus 1ex minus .2ex
                          \skip\afterskip=1.5ex plus .2ex
                          \baselineskip=15pt
        \par
        \def \par{\hss\egroup\line\bgroup\hss}\obeylines
        \line\bgroup\hss
\global\advance\thesubsection by 1
\xdef \lastref{\number\thesection.\number\thesubsection}
\global\thesubsubsection=0 \global\theparagraf=0
\number\thesection.\number\thesubsection\quad #1\hss\egroup\endgroup \par}

\def\\{\par}

\long\def\subsubpartie#1{\bgroup\bf
                                \par\penalty -100
                                \vskip 3.25ex plus 1ex minus .2ex
                                \skip\afterskip=1.5ex plus .2ex
                                
\begincentered
\global\advance\thesubsubsection by 1
\xdef \lastref{\number\thesection.\number\thesubsection.\number\thesubsubsection}
\global\theparagraf=0
\number\thesection.\number\thesubsection.%
\number\thesubsubsection\quad #1\endcentered\egroup\par}
\def\paragrafsubsub{%
\global\advance\theparagraf by 1
\xdef \lastref{\number\thesection.\number\thesubsection.\number\thesubsubsection.%
\romannumeral\theparagraf}
{\bf \number\thesection.\number\thesubsection.\number\thesubsubsection.\romannumeral\theparagraf}%
\kern0.2em --- }

\def\paragrafsub{%
\global\advance\theparagraf by 1
\xdef \lastref{\number\thesection.\number\thesubsection.\number\theparagraf}
{\bf \number\thesection.\number\thesubsection.\number\theparagraf}%
\kern0.2em --- }

\def\paragraf{%
\global\advance\theparagraf by 1 
\xdef \lastref{\number\thesection.\number\theparagraf}
{\bf \number\thesection.\number\theparagraf}%
\kern0.2em --- }

\def\paragraphe{%
\par \indent
\ifcase\thesubsection %
  \paragraf
\else
\ifcase\thesubsubsection\paragrafsub %
 \else\paragrafsubsub\fi
\fi}

\def\numeq{\global\advance\theequ by 1%
        \xdef \lastref{(\number\theequ)}%
        \eqno{(\number\theequ)}}
\def\nume{\global\advance\theequ by 1 
        \xdef \lastref {(\number\theequ)}%
        (\number\theequ)}
\def\theo{\global\advance\thetheo by 1%
\xdef\lastref{\number\thetheo}%
\number\thetheo}

\newwrite\fileref
\newread\instream
\newif \ifrefmodif 
\def \defineref#1#2{{\def\next{#1}%
        \expandafter\xdef
           \csname ref = \meaning\next\endcsname{#2}%
        }}

\def \initfileref{
\openin\instream=\jobname.ref
\ifeof\instream \message{ Ahem } 
                \message{**********************************************************}%
                \message{******* Le fichier \jobname.ref n'existait pas ***********}%
                \message{********** Il faudra absolument recompiler ***************}%
                \message{**********************************************************}%
                \message{ }%
\else \closein\instream \input \jobname.ref \refmodiffalse   
\fi
\ifarxiv\else\immediate\openout\fileref=\jobname.ref\fi%
\global\let\initfileref=\relax} 

\def \label#1{\ifarxiv\initfileref\else{%
        \toks0={#1}\wlog{REF \the\toks0= \lastref}
        \initfileref
        \immediate\write\fileref{\noexpand\defineref%
                {\the\toks0}{\lastref}}%
        \def\next{#1}%
        \expandafter\ifx 
          \csname ref = \meaning\next\endcsname\lastref
        \else \global\refmodiftrue  \message{ }\message{ Attention }
        \message{reference {\the\toks0 = \lastref }  modifiee ou redefinie} \message{ }\fi 
        \defineref{#1}{\lastref}%
        \ifrefvis\ifhmode\raise 6pt \vbox{\hbox to 0pt{\hss\fivebf [\the\toks0]\hss}}
        \else\llap{\fivebf [\the\toks0]}\fi\fi}\fi}

\def \ref#1{{\initfileref \def\next{#1}%
\csname ref = \meaning\next\endcsname}}

\def\cite#1{{\initfileref \def\next{#1}%
\csname ref = \meaning\next\endcsname}}


\def \ignorepar{\afterassignment\ignoreparaux \let\next=}
\def \ignoreparaux{\ifx\next\par \let\next\ignorepar \fi \next}

\def\bibliography{\bgroup\sectionfont
        \par\penalty -500
        \vskip 3.25ex plus 1ex minus .2ex
        \skip\afterskip=1.5ex plus .2ex
\begincentered References \endcentered\egroup\par
\tabskip5pt minus 1pt%
\noindent\halign to \hsize \bgroup##\hfil&##\hfil\cr}
\def\bibliend{\egroup}

\def\bibitem[#1]#2{\xdef\lastref{[#1]}\label{#2}\hskip-8pt[#1]&%
\vtop\bgroup\hsize10.6cm\noindent\ignorespaces}
\def\finitem{\medskip \egroup \cr}

\def\newblock{\hskip .11em plus .33em minus .07em}

\def\em{\sl}

\def\proof{\noindent{\bf Proof. }}
\def\qed{\ \-\hfill {\ding \char 111} \par}

\catcode`\@=11
\def\hexnumber@#1{\ifcase#1 0\or 1\or 2\or 3\or 4\or 5\or 6\or 7\or 8\or
 9\or A\or B\or C\or D\or E\or F\fi}

\font\tenmsa=msam10
\font\sevenmsa=msam7
\font\fivemsa=msam5
\newfam\msafam
\textfont\msafam=\tenmsa
\scriptfont\msafam=\sevenmsa
\scriptscriptfont\msafam=\fivemsa
\edef\msafam@{\hexnumber@\msafam}
\mathchardef\dabar@"0\msafam@39
\def\dashrightarrow{\mathrel{\dabar@\dabar@\mathchar"0\msafam@4B}}
\def\dashleftarrow{\mathrel{\mathchar"0\msafam@4C\dabar@\dabar@}}

\def\ulcorner{\delimiter"4\msafam@70\msafam@70 }
\def\urcorner{\delimiter"5\msafam@71\msafam@71 }
\def\llcorner{\delimiter"4\msafam@78\msafam@78 }
\def\lrcorner{\delimiter"5\msafam@79\msafam@79 }
\def\yen{{\mathhexbox@\msafam@55}}
\def\checkmark{{\mathhexbox@\msafam@58}}
\def\circledR{{\mathhexbox@\msafam@72}}
\def\maltese{{\mathhexbox@\msafam@7A}}

\catcode`\@=12

\def\T{{\Bbb T}}

\def\R{{\Bbb R}}

\def\frac#1#2{{#1\over#2}}

\def \prodscal#1,#2>{\langle #1,#2 \rangle}

%
%




\arXiv 

\anglais
\tenpoint

\TITR The Length of Harmonic Forms on a Compact Riemannian Manifold|


\auteur P.A. Nagy and C. Vernicos*|

\vfootnote*{\sevenrm The second author was partially supported by european project ACR OFES number 00.0349 and grant of the FNRS 20-65060.01}



\centerline{\bf Abstract}
\medskip
\centerline{\vbox{\hsize 10cm \baselineskip=2.5ex \hautdepage We study $n$ dimensional
Riemanniann
manifolds with harmonic forms of constant length and first Betti number
equal to $n-1$ showing that they are $2$-steps nilmanifolds with some special metrics. 
We also characterise, in terms of properties on the product of harmonic forms, the
left invariant metrics among them.
This allows us to clarify the case of equality in the stable isosytolic inequalities
in that setting.
We also discuss other values of the Betti number.
}}

\vfootnote{}{{\sevensl Key-words : }{\ssubhtdpfont Harmonic forms, Spectrum of
the Laplacian.}}
  
\vfootnote{}{{\sevensl Classification : }{\ssubhtdpfont 53C20, 58J50 }}

\partie{Introduction}

Let $(M^n,g)$ be a Riemannian manifold. In this note we focus on the case where
all harmonic forms are of constant length. Recently these manifold have appeared
in different settings where they play a singular part.

In dimension $4$, recent work of Lebrun \cite{lebrun} shows a strong interplay
between the length of harmonic self-dual $2$-forms of the manifold and the non-vanishing
of Seiberg-Witten invariants, in particular the existence of a symplectic structure.

{\bf Geometrically Formal Manifolds}, are closed Riemaniann manifolds having a metric
such that the space of harmonic forms is a sub-algebra of the algebra of
differential forms. If the manifold is
also oriented, one easily sees that all harmonic forms have constant pointwise
norm (see \cite{kot}). In particular if the first Betti number is equal to the dimension (which
actually bounds it), then the manifold is a flat torus (In fact it
is also true if some other Betti number is maximal, provided the dimension is prime
for example). 
From a remark made by Kotschik \cite{kot} we also know that being geometrically
formal and orientable implies a second obstruction on the first Betti number, 
if $n$ is the dimension,
the Betti number can't be $n-1$. However, $n$ dimensional manifolds
with one-harmonic forms of constant length, and first Betti number equal to $n-1$ exists. 

The main result of this paper is the following theorem which says exactly what these manifolds
are (see Definition \ref{defone} for the meaning of pseudo left invariant) 

\proclaim Theorem \theo. Let $(N^{n+1},g)$ be a compact orientable connected manifold such
that all of its harmonic $1$-forms are of constant length, and
such that $b_1(N^{n+1})=n$, then $(N^{n+1},g)$ is a $2$-step nilmanifold
whose kernel is one dimensional and $g$ is a pseudo left invariant metric.\label{theonil1}

{\bf The macroscopical spectrum of a nilmanifold}, is given by the asymptotic behaviour
of the eigenvalues of Laplace-Beltrami operator acting on the function on the metric balls
of the universal covering of a nilmanifold, as the radius of the balls goes to infinity.
In \cite{vernicos} the second author showed that the first eigenvalue of this
macroscopical spectrum satisfies an inequality, whose equality case is attained by the nilmanifolds
having all harmonic $1$-forms of constant norm. This shows that the nilmanifolds with left-invariant
metrics are not the only one satisfying the equality case, as in the torus case. 
Hence theorem \ref{theonil1} gives the following corollary in that setting :

\proclaim Corollary \theo. Let $(M^n,g)$ be a nilpotent Riemannian manifold, with first
Betti number $b_1=n-1$. Let $B_g(\rho)$ be the ball of radius $\rho$ induced by the lifted metrics
on the universal covering of $M^n$. Let $\lambda_1\bigl(B_g(\rho)\bigr)$ be the first eigenvalue of the laplacian
acting on functions over $B_g(\rho)$ for the Dirichlet problem. Then there are some functions
$\lambda_1^\infty(g)$ and $\lambda_1^{al}(g)$ such that
$$
\lim_{\rho \to \infty} \rho^2\lambda_1\bigl(B_g(\rho)\bigr) = \lambda_1^\infty(g) \leq \lambda_1^{al}(g)
$$
with equality if and only if $M^n$ is a $2$-step nilmanifold with one dimensional center 
and $g$ is pseudo left invariant.

In this
paper we also study what other assumption, in term of product of harmonic forms (lead
by the geometrically formal background), one could add to characterize the left invariant metrics
in the case of the Betti number equal $n-1$, where $n$ is the dimension (see theorem \ref{theonil2}
to \ref{theonil4}).

{\bf Isosystolic stable inequalities}, studied among other by V.~Bangert and M.~Katz \cite{bk}
(see also references therein), give lower bounds on the volume of compact orientable manifold
in terms of some short closed geodesics (systoles). The cases of equality is obtained by manifolds
with one-forms of constant length. Hence Theorem \ref{theonil1} implies the following corollary
(where ${\rm stsys}_1$ is the stable systole) 

\proclaim Corollary \theo. Let $(X,g)$ be a compact, oriented, $n$-dimensional Riemannian
manifold with first Betti number $b=n-1$. Then there is a constant $c(n)$ such that
$$
{\rm stsys}_1(g){\rm sys}_{n-1}(g)\leq c(n) {\rm vol}_g(X)
$$
Equality occurs if and only if $X$ is a $2$-step nilmanifolds and there exist
a dual-critical lattices $L$ in $\R^{n-1}$ and a submersion of $X$ into $\R^{n-1}/L$ (with
minimal fibers).

Hence in the three dimensional case, equality is only satisfied by the three dimensional
Heisenberg group endowed with a metric such that there is 
a Riemannian submersion onto an equilateral torus.

Our result are also to be compared with the work of E.~Aubry, B.~Colbois, P.~Ghanaat and
E.A. Ruh \cite{acgr}, where it is showed that an $n$-dimensional oriented manifold $M$ having
$n-1$ small (compared to the diameter) eigenvalues 
(for the laplacian acting on one-forms) is diffeomorphic
to a nilmanifold with an almost left invariant metric. However, in their paper the
authors needed strong assumptions on the curvature. In this note, instead of a control
of the curvature, we have a control of the length of harmonic forms. We would like
to stress that this seems to be a hidden assumption on the curvature.

\partie{Implication of the existence of harmonic forms

of constant length.}

\subpartie{The Albanese map}

Let $(M^n,g)$ be a compact Riemannian manifold having all of its harmonic one forms
of constant length. 

Let us stress that we do not consider the case where there
are harmonic forms of constant length, but that {\sl all} harmonic forms
are supposed to be of constant length. This implies for example that the pointwise
scalar product of two harmonic forms is constant, which is not the case with just
the existence assumption as the following example shows (which happens to answers
question $7$ of section $10$ in \cite{bk}:

\noindent{\bf Example \theo.} Let $(M^n,g)$ be  Riemannian manifold, consider
$N^{n+2}=M^n\times\T^2$ endowed with the following Riemannian metric $h=g\oplus s$
where $s$ is defined as follows. Take a function from $M$ to $\R$ such
that $f^2< 1$ then we take the metric $s$ on $\T^2$ such that $s(dx,dx)=1=s(dy,dy)$ and
$s(dx,dy)=f$.
We claim that $\alpha_1=dx$ and $\alpha_2=dy$  which are closed are also  co-closed thus harmonic,
but by construction their scalar product is not constant.

\proof
with respect to $h$ in $N^{n+2}$. Let us take  $(e_i)_{1\leq i\leq n+2}$  a local orthonormal basis
of $N^{n+2}$, with $(e_k)_{1\leq k\leq n}$ spanning $TM$, and $e_{n+1}$, $e_{n+2}$ spanning $T\T^2$. 
Consider also $(X_i)_{i=1,2}$ the dual vector field with respect to $h$ of $(\alpha_i)_{i=1,2}$.
Then remark that $X_i = x_i^1\frac{\partial}{\partial x} + x_i^2\frac{\partial}{\partial y}$ and
$e_k=E_k^1\frac{\partial}{\partial x} + E^2_k\frac{\partial}{\partial y}$ 
where $(x_i^j)_{i,j=1,2}$ and $(E_k^i)$ are functions from $M$ to $\R$.
Now write
$$
-d^*\alpha_i=\sum_{1\leq k\leq n+2} h(\nabla_{e_k}X_i,e_k)
$$
As the Levi-Civita connection is torsion free, we remark that for $i=1,2$ and $1\leq k\leq n+2$
$$
h(\nabla_{e_k}X_i,e_k) =  h(\nabla_{X_i}e_k,e_k) +h([X_i,e_k],e_k) {\rm .}
$$
However for $i=1,2$ and $1\leq k\leq n+2$
$$
h(\nabla_{X_i}e_k,e_k) = \frac{1}{2} X_i\cdot h(e_k,e_k) =0
$$
and noticing that for $1\leq k\leq n$, $[\frac{\partial}{\partial x},e_k]=[\frac{\partial}{\partial y},e_k]=0$ we obtain
$$\eqalign{
[X_i,e_k]&= [x_i^1\frac{\partial}{\partial x} + x_i^2\frac{\partial}{\partial y}, e_k ]\cr 
      &=-(e_k\cdot x_i^1) \frac{\partial}{\partial x} -  (e_k\cdot x_i^2) \frac{\partial}{\partial y}}
$$
and has $\frac{\partial}{\partial x}$ and $\frac{\partial}{\partial y}$ are in $T\T^2$ we get that for all $1\leq k\leq n$ and $i=1,2$,
$$
h(\nabla_{e_k}X_i,e_k) = 0
$$
Now for $k=n+1,n+2$ we have
$$
[X_i,e_k]=[x_i^1\frac{\partial}{\partial x} + x_i^2\frac{\partial}{\partial y}, E_k^1\frac{\partial}{\partial x} + E^2_k\frac{\partial}{\partial y}]
$$
however, for any function $f$ define on $M$, $T\T^2 \in \ker df$ and noticing
that $[\frac{\partial}{\partial x},\frac{\partial}{\partial y}]=0$ we finally get that 
$\alpha_i$ is indeed harmonic.
\qed

Let $b_1$ be the first
Betti number of our manifold $M$. Using a basis of harmonic $1$-forms 
we can define the Albanese map (or Jacobi map) $\pi$ using integration,
this gives a map onto a torus $\T^{b_1}$, on which we put the usual flat metric. 
Besides
we have $b_1\leq n$, because we can't have more than $n$ linearly independant
$1$-forms. Hence in this case $\pi$ is a Riemannian submersion.
Moreover Albanese's map being harmonic (see A.~Lichnerowicz \ref{lich}), 
the fibers are minimal (see J.~Eells and J.M. Sampson \ref{eesa} and 9.34 \ref{besse} page 243).

A theorem of R.~Hermann (see theorem 9.3 in \ref{besse} page 235) states
that a Riemannian submersion whose total space is complete is a locally trivial
fibre bundle. We sum this up in the following proposition.

\proclaim Proposition \theo. Let $(M^n,g)$ be a compact Riemannian manifold and $b_1$ be its
first Betti number. Then all harmonic $1$-forms are of constant length if and only if $(M^n,g)$ is
a locally trivial fiber bundle, with minimal fibers, over a $b_1$-dimensional flat torus, $b_1\leq n$ 
$$
F^{n-b_1} \hookrightarrow M^n {\buildrel \pi\over \longrightarrow} \T^{b_1}
$$
Moreover if $b_1=n$ then $\pi$ is a Riemanian isometry,
hence $(M^n,g)$ is a flat torus.\label{prop1}






Now let us look at some other consequences of
the existence of harmonic forms of constant length. Thanks to the Albanese map
we can lift the harmonic forms of the Albanese torus on the manifold. Let us call
$\alpha_1,\dots,\alpha_{b_1}$ an orthonormal family of lifted harmonic forms. Using the duality through the metric
we can associates to these harmonic forms $b_1$ vector fields $X_1,\dots, X_{b_1}$.
These vector fields define a sub-bundle ${\cal H}$, which we will call the {\bf Horizontal},
such that for all $x\in M$, ${\cal H}_x$ is  generated by $X_1(x),\dots, X_{b_1}(x)$.

If we let ${\cal V}$ be the orthogonal complement of ${\cal H}$ with respect to the metric, 
which we will refer to as the {\bf Vertical}, then we have the following 

\proclaim Proposition \theo. Let $(M^n,g)$ be a compact Riemannian manifold and $b_1$ be its
first Betti number. If all harmonic $1$-forms are of constant length, then there is
distribution ${\cal H}$ given by an orthonormal frame of vector fields $X_1,\dots,X_{b_1}$ 
dual to an orthonormal frame of
harmonic $1$-forms, such that the tangent bundle splits orthogonaly as follows :
$$
TM={\cal V}\oplus{\cal H} { \rm.}
$$
Moreover the distribution ${\cal V}$ is integrable, and for any $1\leq i,j\leq b_1$ and $U\in {\cal V}$, 
$[X_i,X_j] \in {\cal V}$ and $[X_i,U]\in {\cal V}$.\label{hetv}

\proof
This comes from the fact that the forms $(\alpha_i)_{1\leq i\leq b_1}$ are closed. Indeed for any closed $1$-form
$\alpha$ we have the following equality for any $X$, $Y$ in $TM$ :
$$
\alpha([X,Y]) = X\cdot\alpha(Y) - Y\cdot\alpha(X)
$$
thus for any $i$, $j$ and $k$ we have :
$$
\alpha_i([X_j,X_k])= X_j\cdot\alpha_i(X_k) -X_k\cdot\alpha_i(X_j) = X_j\cdot\delta_{ik} -X_k\cdot\delta_{ij} =0
$$
hence $[X_j,X_k]$ is orthogonal to any $X_i$. If $U$ and $V$ are vertical vector fields,
then it is easily seen that for any $1\leq i\leq b_1$
$$
\alpha_i([X_j,U]) = 0, \quad {\rm and} \quad  \alpha_i([U,V])= 0
$$
\qed

\subpartie{An useful decomposition}
Let $(M^n,g)$ be a compact Riemannian manifold with a unit vector field $Z$. Let ${\cal V}$
be the distribution generated by $Z$ (which is sometimes called the {\bf Vertical} distribution)
and ${\cal H}$ its orthogonal complement with respect to $g$ (wich we will call {\bf Horizontal}).
Then the tangent bundle splits as follows
$$
TM={\cal V}\oplus {\cal H}.
$$ 

If $i_Z(\cdot)=Z\lrcorner\cdot$, is the interior product by $Z$ then we can define
the space of horizontal $p$-forms as follows :
$$
\Lambda^p({\cal H})=\Lambda^p(M)\cap\ker(i_Z).
$$

Let us introduce the $1$-form $\vartheta=Z^\flat$ dual to $Z$ with respect to $g$. Then it is
an easy exercise to see that we have the following decomposition
of $p$-forms :
$$
\Lambda^p(M)=\Lambda^p({\cal H})\oplus\bigl[\Lambda^{p-1}({\cal H})\land\vartheta\bigr]. \numeq \label{decomposition}
$$ 

Notice that following this decomposition we have 
$$
d\vartheta=b + \eta\land\vartheta
$$
where $b$ is a horizontal $2$-form called the {\bf curvature} of the horizontal distribution $H$
and $\eta$ is the horizontal one form associated to the horizontal vector field
$-\nabla_ZZ$ with respect to $g$, a fact the reader may easily verify.

Let us introduce the horizontal exterior differential $d_H$, which
associates to an horizontal differential form the horizontal part of
its exterior differential.

We also introduce the multiplication operators
$$\eqalign{%
L &: \Lambda^q({\cal H}) \to \Lambda^{q+2}({\cal H}), \quad L:=\cdot\land b \cr
S &: \Lambda^q({\cal H}) \to \Lambda^{q+1}({\cal H}), \quad S:=\eta\land \cdot } \numeq
$$

Thanks to the decomposition \ref{decomposition}, each $p$-form may be identified to
a couple $(\alpha,\beta)\in \Lambda^p({\cal H})\times\Lambda^{p-1}({\cal H})$ of horizontal forms. Thus we can see the exterior 
differential $d$ acting on $p$-forms as an operator from $\Lambda^p({\cal H})\times\Lambda^{p-1}({\cal H})$ 
to $\Lambda^{p+1}({\cal H})\times\Lambda^{p+2}({\cal H})$. Then, we have the
following proposition, where ${\cal L}_Z$ is the Lie derivative in the direction of $Z$,
and where for any differential operator $D$, $D^*$ is its formal adjoint with respect to $g$.

\proclaim Proposition {\theo.} \cite{nagy}.  With respect to the decomposition \ref{decomposition}
we have for the exterior differential acting on $p$-forms : \label{lastprop}
$$
d = \pmatrix{ d_H & (-1)^{p-1}L \cr (-1)^p{\cal L}_Z & d_H+S} 
$$
and for the codifferential :
$$
d^* =\pmatrix{ d^*_H & (-1)^{p-1}{\cal L}_Z^* \cr (-1)^pL^* & d_H^*+S^* }
$$




\partie{The case of the first Betti number 

being one less that the dimension}

The aim of this section is to describe the $n$-dimensional manifolds having all of their
harmonic one-forms of constant length, and with first Betti number equal to $n-1$.
All manifolds considered here are of dimension greater or equal to $3$ (in
dimension $1$ and $2$ we have only the circle $S^1$ and the torus $\T^2$ which admit
forms of constant length). 

Let us recall that all $2$-step nilmanifolds are principal torus bundles over a torus
(see R.S.~Palais and T.E.~Stewart \cite{palais})

\proclaim Definition \theo. Let $N^{n+1}$ be a $2$-step nilmanifold whose center
is one dimensional and
consider a submersion $p$ from $N^{n+1}$
onto a flat torus $\T^n$. Let $(\alpha_1,\dots,\alpha_n)$ be the lift of an orthonormal basis of harmonic
$1$-forms over the torus. 
Choose a principal connection form $\vartheta$ for this submersion. 
Let $g_\vartheta$ be the Riemanniann metric
such that the dual basis of $(\alpha_1,\dots,\alpha_n,\vartheta)$ is orthonormal. Thus
$p$ becomes a Riemannian submersion. We will call such a metric {\rm pseudo left invariant}.%
\label{defone}

Remark that if the dual basis of $(\alpha_1,\dots,\alpha_n,\vartheta)$ in this definition is left invariant,
then the metric $g_\vartheta$ is left invariant.

Thanks to this definition we can recall our characterization

\proclaim Theorem \ref{theonil1}. Let $(N^{n+1},g)$ be a compact orientable connected manifold such
that all of its harmonic $1$-forms are of constant length, and
such that $b_1(N^{n+1})=n$, then $(N^{n+1},g)$ is a $2$-step nilmanifold
whose kernel is one dimensional and $g$ is a pseudo left invariant metric.\label{theonil1}

\proof
We deduce from proposition \ref{prop1} that we have the following fibration :
$$
S^1 \hookrightarrow N^{n+1} {\buildrel \pi\over \longrightarrow} \T^n
$$ 
where $\pi$ is the Albanese map.


We begin by taking a basis $(a_1,\dots,a_n)$ of harmonic forms over the Albanese torus
and lift it  to a basis of harmonic $1$-forms $(\alpha_1,\dots,\alpha_n)$ over $N^{n+1}$. Let $X_1,\dots,X_n$
be their dual vector fields with respect to the metric, they span, by definition, ${\cal H}$.

Thanks to Proposition \ref{hetv}, for any $i,j$ 
$$
[X_i,X_j] \in {\cal V} \numeq \label{crochet1} 
$$
(remark that we can also deduce from that fact that $[X_i,X_j]=2\nabla_{X_i}X_j$).
Now take $Z$ the dual vector field to the $1$-form $Z^\flat=*(\alpha_1\land\dots\land\alpha_n)$ ($*$ is the Hodge
operator, thus this form is co-closed), its length is
constant by construction. Furthermore $Z$ belongs to and spans ${\cal V}$.

For any co-closed one form $\alpha$ we have
$$
\sum_k \bigl(\nabla_{X_k}\alpha\bigr) X_k + \bigl(\nabla_Z\alpha\bigr) Z =0 \numeq \label{coclcond}
$$
Taking for $\alpha$ each of the $\alpha_i$ in turn, since $[X_i,X_j]=2\nabla_{X_i}X_j\in{\cal V} $, 
we deduce from equality \ref{coclcond} that for $i=1,\dots, n$
$$
g(\nabla_ZX_i, Z) =0
$$
however as the Levi-Civita is torsion free, for $i=1,\dots,n$
$$g(\nabla_ZX_i, Z) = g(\nabla_{X_i}Z, Z) + g([Z,X_i],Z) = g([Z,X_i],Z) $$
but by proposition \ref{hetv}, $[Z,X_i] \in {\cal V} $ hence for any $i=1,\dots,n$
$$
[X_i,Z] = 0 \numeq \label{crochet2}
$$
This also implies that $Z$ is a Killing field.

From \ref{crochet1} we have the existence of functions $f_{ij}$ such that
$$
[X_i,X_j] = f_{ij} Z {\rm .} \numeq \label{crochet3}
$$

However we would like to have some structural constants instead of
the functions $f_{ij}$. 

Let us remark that ($Z$ being a Killing field) we have
$$
dZ^\flat(X,Y)=2g(\nabla_XZ,Y) \numeq \label{crochet4}
$$
thus if we decompose $dZ^\flat$ in the basis given by $\alpha_i\land\alpha_j$ and $Z^\flat\land\alpha_i$  for all $i,j$ then
thanks to \ref{crochet3} and \ref{crochet4} we get that
$$
dZ^\flat=\sum_{i<j} f_{ij} \alpha_i\land\alpha_j \numeq \label{pasconstantes}
$$
In other words $dZ^\flat$ is horizontal and because $Z$ is a Killing field it is projectable, i.e. , 
there exists a unique $2$-form $\beta$ on the Albanese torus such that
$dZ^\flat=\pi^*\beta$.

Remark that $d\beta=0$, thus $\beta=\beta_0+d\alpha$ by the Hodge-de Rham theorem, with $\alpha\in \Lambda^1(\T^{b_1})$
and $\beta_0$ harmonic.
Hence if $\zeta_0=Z^\flat-\pi^*\alpha$, then $d\zeta_0=\pi^*\beta_0$. But now
$\beta_0=\sum_{ij} c_{ij} a_i\land a_j$ where the $c_{ij}$ are constants.
This implies that
$$
d\zeta_0=\sum_{ij} c_{ij} \alpha_i\land\alpha_j \numeq \label{enfinconstant}
$$

Notice that $\zeta_0$ is also a (principal) connection $1$-form.
We are now taking as a basis of vector fields the dual base $(X_1^0,\dots,X_{n-1}^0,Z_0)$ of
the base $(\alpha_1,\dots,\alpha_{n-1},\zeta_0)$ (i.e. $\alpha_i(X_j^0)=\delta_{ij}$, $\ker \zeta_0=\langle X_1^0,\dots,X_{n-1}^0\rangle$ and $\zeta_0(Z_0)=1$)  
then we have from \ref{enfinconstant}
$$\eqalign{
[X_i^0,Z_0]&=0\cr
[X_i^0,X_j^0]&= c_{ij} Z_0} 
$$

Thus we can built an homomorphism $I_A$ between the Lie algebra $\mathfrak{a}$ 
defined by $A_i,\dots,A_n, A_{n+1}$ and with brackets
$$
[A_i,A_j] =c_{ij} A_{n+1}
$$  
all the other brackets being equal to zero, by taking $I_A(A_i)=X_i^0$ for $i=1,\dots,n$
and $I_A(A_{n+1})=Z_0$.

Now if ${\cal A}$ is the simply connected Lie group associated to $\mathfrak{a}$,
thanks to the compactness of $N^{n+1}$ we can integrate the homomorphisme $I_A$ to
obtain an action of ${\cal A}$ on $N^{n+1}$ (see Corollary $3$ and $4$ of
theorem 2.9 page 113 in \ref{ency}) . 

Since each orbit is open and $N^{n+1}$ is
connected, this action is transitive. From this we deduce 
(see F.W.~Warner \ref{warner} theorem 3.62) that $N^n$ is a Lie group. And thanks
to the constant of structures $(c_{ij})$ we deduce that it is a two step nilpotent lie
group (see M.~Spivak \cite{spivak} volume I theorem 17 in chapter 10 for example).

Concerning $g$, it is pseudo left invariant by the above discussion.





\qed

It is worth noticing that if the manifold is not orientable this is no more
true. Indeed, in the three dimensional case there are quotients of flat tori
by a finite group whith first Betti number equal to $2$.

Note also that theorem \ref{theonil1} implies
that if $M^n$ is a geometrically formal closed oriented manifold
then $b_1(M^n)\not= n-1$ as already showed by D.~Kotschik \cite{kot}.

A natural question arising is whether one can characterize
the left invariant metrics among pseudo left invariant metric ones,
in term of properties of harmonic forms.
We first look at the $3$-dimensional case.

\proclaim Theorem \theo. Let $(N^3,g)$ be a compact orientable connected manifold such
that all harmonic $1$-forms are of constant length, $b_1(N^3)=2$, and such that
the wedge product of two harmonic $1$-forms
is an eigenform of the laplacian, then $(N^3,g)$ is a compact quotient of the $3$-dimensional
Heisenberg group and $g$ is a left invariant metric. \label{theonil2}

\proof
From the added assumption we get that $Z^\flat$ of the previous
proof is an eigenform of the laplacian.
That is to say that there is some constant $\lambda$ such that
$$
\Delta Z^\flat= \lambda Z^\flat
$$
However, using the useful decomposition, knowing that $Z^\flat$ is coclosed
and $dZ^b=b=\pi^*(\beta)$ (following the decomposition \ref{decomposition}, $b$
corresponds to $(b,0)$ and $Z^\flat$ corresponds to $(0,1)$ ) we must also have that
$$
\Delta Z^\flat=d^*dZ^b= d^*b = (L^*b) Z^\flat  = (L^*L\cdot 1) Z^\flat = |\beta|^2 Z^\flat 
$$
which means that $\beta$ is of constant length. But $\beta$ is a $2$-form on
a $2$-dimensional Torus, which also means that $\beta$ is proportional
to the volume form, in other words the function (which is unique in that case)
 in the equality \ref{pasconstantes}
is a constant. Remark that this also tells us that the eigenvalue is $\lambda=|\beta|^2=f_{12}^2$.
\qed

Now for the higher dimensional case.

\proclaim Theorem \theo. For $n>3$, let $(N^n,g)$ be a compact orientable connected manifold such
that all harmonic $1$-forms are of constant length, $b_1(N^n)=n-1$
and the wedge product of any $n-2$ harmonic 
$1$-forms is an eigenform of the Laplacian,
then $(N^n,g)$ is a $2$-step nilmanifold
whose kernel is one dimensional and $g$ is a left invariant metric. \label{theonil3}

\proof We use the same notations as in the proof of theorem \ref{theonil1}.
From the new assumption we get thanks to  the star Hodge operator
that for any harmonic $1$-form $\alpha$, $\alpha\land Z^\flat$ is a co-closed eigenform. Using
the decomposition \ref{decomposition} we associate to $\alpha\land Z^\flat$ the pair
$(0,\alpha)$. As $\alpha\land Z^\flat$ is coclosed $\Delta(\alpha\land Z^\flat)=d^*d(\alpha\land Z^\flat)$. Let us use
proposition \ref{lastprop} without forgetting that $\alpha$ is harmonic, and noticing
that as $dZ^\flat$ is horizontal by the proof of theorem \ref{theonil1}, $S=0$
$$
d(\alpha\land Z^\flat)\equiv\pmatrix{ d_H & -L \cr {\cal L}_Z & d_H} \pmatrix{0\cr \alpha} = \pmatrix{-L \alpha\cr d_H\alpha}%
=\pmatrix{-L(\alpha)\cr 0}
$$
and 
$$
d^*d(\alpha\land Z^\flat)\equiv\pmatrix{ d^*_H & {\cal L}_Z^* \cr -L^* & d_H^* }\pmatrix{-L(\alpha)\cr 0} = 
\pmatrix{d_H^*\bigl(-L(\alpha)\bigr) \cr L^*L(\alpha)}
$$
but the last term is also equal to $(0,\lambda\alpha)$ by the eigenform assumption,
hence
$$
d_H^*(\alpha\land dZ^\flat)=0
$$

But $\alpha$ and $dZ^\flat$ are horizontal and projectable, thus in fact there are
one $1$-form $a$ and one $2$-form $\beta$ such that we can write
$$
0= d_H^*(\alpha\land dZ^\flat) = \pi^*\bigl(d^*(a\land\beta)\bigr).
$$

In other words, for any $1$-harmonic form $a$ on the torus we have
$$
d^*(a\land\beta)=0. \numeq \label{awbcoferme}
$$
Let us take an orthonormal base $(e_i)$ of parallel vector fields, then for any form $\omega$ 
$$
-d^*\omega= \sum_i e_i\lrcorner\nabla_{e_i}\omega 
$$
We want to apply that last equality to $a\land\beta$. First notice that ($a$ is parallel)
$$
\nabla_{e_i} (a\land\beta) = a\land\nabla_{e_i}\beta.
$$

Contracting  by $e_i$ one obtains
$$
e_i\lrcorner\bigl(\nabla_{e_i} (a\land\beta)\bigr) = a(e_i)\nabla_{e_i}\beta - \alpha\land(e_i\lrcorner\nabla_{e_i}\beta)
$$
Sum over $i$ and use the coclosed condition (i.e. \ref{awbcoferme})
and you get
$$
\nabla_{a^\#}\beta + a\land d^*\beta =0 \numeq \label{reduction1}
$$
Now we take the interior product with $e_i$ of \ref{reduction1}, with $a=e_i^\flat$ which gives
$$
e_i\lrcorner\nabla_{e_i}\beta + e_i\lrcorner(e_i^\flat\land d^*\beta)=
e_i\lrcorner\nabla_{e_i}\beta + d^*\beta - d^*\beta(e_i)e_i^\flat =0
$$
we sum over $i$ one more time
$$
(n-3) d^*\beta=0
$$
and by assumption $n>3$ thus $d^*\beta=0$, that is to say that $\beta$ is harmonic
over the torus. Hence it follows that $\beta$ is parallel 
and all the functions $f_{ij}$ of \ref{pasconstantes} are constants,
which allows us to conclude. 
\qed

There is another case where things can be precised.

\proclaim Theorem \theo. Let $(N^{2m+1},\omega,g_\omega)$ be a compact contact manifold with a contact
form $\omega$ and an adapted
Riemannian metric $g_\omega$ such that all harmonic $1$-forms are of constant length,
then $(N^{2m+1},g)$ is a compact quotient of an Heisenberg group and $g$ is a left invariant metric.%
\label{theonil4}

\proof
From theorem \ref{theonil1} we get that $N^{2m+1}$ is a two step nilmanifold
and $g_\omega$ is pseudo left-invariant. We also get that $g_\omega=\vartheta^2+\pi^*(h)$ where $\vartheta$ is
a one form, such that $d\vartheta=\pi^*(\beta)$ for some closed $2$-form $\beta$ over the Albanese torus.
Now theorem 3.2 and theorem 3.4 of H.P. Pak and T. Takahashi in \cite{paktak} implies that
for all harmonic $1$-forms $\alpha$, if $T$ is the Reeb vector field attached to $\omega$,
$$
T\lrcorner \alpha = \alpha(T) =0
$$
wich means that $T=fZ$ for some function $f$, and as $T$ and $Z$ are of constant
unit length for the metric $g_\omega$ it means that $\omega=\vartheta$.
Now, thanks to theorem  \ref{theonil1} we know that $T=Z$ is a Killing field.
This implies that the almost complex structure $J$ on $\ker \omega$ lives on
the flat torus given by the Albanese submersion.
Hence we are in front of an almost-K\"ahler flat torus, but following
\cite{olszak} and \cite{armstrong}  it has to be K\"ahler, thus $\beta$ is parallel.
\qed

\partie{Some remarks on the general case}

The aim of this section is to point out the main differences between the case
$b_1=n-1$ and $b_1\leq n-2$ for $n$ dimensional manifolds admitting one-harmonic forms
of constant length. We want to give some hint on the failure of our approach.

Our first remark is that one should restrict oneself to the study of
locally irreducible orientable Riemannian manifolds to avoid the following cases :
the direct product of a sphere of dimension $p>1$ and
a flat torus of dimension $n-p$, gives a manifold whose first Betti
number is $n-p$, whose dimension is $n$ and with $n-p$ harmonic $1$-forms
of constant length.

The second remark is in the following lemma, which shows the limitation of our
method. Indeed to apply the same ideas one needs far stronger assumptions.

\proclaim Lemma \theo. Let $(N^n,g)$ be a compact locally irreducible manifold
such that all harmonic $1$-forms are of constant length,
$b_1(N^n)=n-p$ and possessing a pointwise orthonormal base $(\vartheta_i)$
of the orthogonal complement of the harmonic $1$-forms. Moreover assume that
$(d\vartheta_i)$ are lifts of closed $2$-forms on the Albanese torus. 
Then $(N^n,g)$ is a two step nilpotent nilmanifolds whose kernel
is $n-p$ dimensional.

\proof
$TM={\cal V}+{\cal H}$
where ${\cal H}$ is spanned by $X_1,\dots,X_{n-p}$ the dual vector fields 
to $\alpha_1,\dots,\alpha_{n-p}$, which are lifts of harmonic $1$-forms on the Albanese torus,
and ${\cal V}$ is the orthogonal complement. We associate, thanks to the metric,
the dual vector fields $(Z_k)$ to the $1$-forms $(\vartheta_k)$.
As the $\alpha_i$ are closed we get that for $1\leq i,j\leq n-p$,
$$
[X_i,X_j]\in {\cal V}.
$$
and our assumptions imply that
$$
d\vartheta_k \in \Lambda^2{\cal H}
$$
where $\Lambda^2{\cal H}= \Lambda^2M \cap \bigcap_k \ker i_{Z_k}$ and that
$$
d\vartheta_k = \pi^*(\beta_k)
$$
where $\beta_k$ is a $2$-form on Albanese's torus. However,
$d\beta_k=0$, hence for some harmonic $2$-form $\beta_k^0$ and some $1$-form $a_k$ over
the Albanese torus we have
$$
\beta_k= \beta^0_k + da_k 
$$
Notice, that $\beta^0_k$ is non zero otherwise $\vartheta_k$ would be horizontal, which
is not the case by assumption. Hence if we consider the independent forms (as
one can easily verify)
$$
\theta_k = \vartheta_k-\pi^*(a_k)
$$
then 
$$
d\theta_k=\pi^*(\beta_k^0){\rm .}
$$
hence we can conclude as in the proof of Theorem \ref{theonil1}.
\qed

As the last section involved nilmanifold, we could decide to
focus on a family of compact manifolds close to them, say solvmanifold.
But even that assumption is not enough to clarify the situation, as the next
lemma with the following remark point out.

\proclaim Lemma \theo. Let $(M^n,g)$ be a solvmanifold
all of whose $1$-forms are of constant length and whose first Betti number is $n-2$, 
then we have the following fibration with minimal fibers :\label{bnmoinsdeux}
$$
\T^2 \hookrightarrow M^n {\buildrel \pi\over \longrightarrow} \T^{n-2}
$$

\proof
This comes from the fact that the fibration in Proposition \ref{prop1}
$$
F^2 \hookrightarrow M^n {\buildrel \pi\over \longrightarrow} \T^{n-2}
$$
gives a long exact sequence on the homotopy groups : 
$$
\dots \to \pi_n(F^2)\to \pi_n(M^n) \to \pi_n(\T^{n-2}) \to \pi_{n-1}(F^2)\to \pi_{n-1}(M^n)\to \dots 
$$
Now, $\T^{n-2}$ and $M^n$ have only their fundamental group which is not trivial,
hence we have the following exact sequence on the fundamental groups
$$
0 \to \pi_1(F^2)\to \pi_1(M^n) \to \pi_1(\T^{n-2}) \to 0
$$
and $\pi_k(F^2)$ is trivial if $k>1$. Which means that $\pi_1(F^2)$ can be seen as a subgroup
of the solvable group $\pi_1(M^n)$, hence it is also solvable.
Thus the fiber is compact
and with solvable fundamental group. However in dimension $2$ the only compact oriented manifold with
solvable fundamental groups are the sphere and the torus, but here the sphere is excluded
because $\pi_2(S^2)\not=0$.
\qed

Without further assumptions we can't expect a more precise results. Indeed the example
of the Sol geometry in dimension $3$, or of any $2$-step nilmanifold with a $2$ dimensional
center will satisfy  Lemma \ref{bnmoinsdeux} with many different metrics,
however built in the same way : following the construction $8.1$ in \cite{bk}.


As a conclusion, the rigidity of  the cases $b_1=n$ and $b_1=n-1$ do not propagates
to lower values of the first Betti number

\bibliography

\bibitem[Arm02]{armstrong}
\bgroup\bf J.~Armstrong\egroup{}.
\newblock An ansatz for almost-K\"ahler, Einstein 4-manifolds.
\newblock {\em J. reine angew. Math.}, 542:53--84, 2002.
\finitem

\bibitem[ACGR01]{acgr}
\bgroup\bf E.~Aubry\egroup{}, \bgroup\bf B.~Colbois\egroup{}, \bgroup\bf
  P.~Ghanaat\egroup{}, and \bgroup\bf E.A. Ruh\egroup{}.
\newblock Curvature, Harnack's inequality, and a spectral characterization of
  nilmanifolds.
\newblock {\em to appear in Ann. of Global Analysis and Geometry}.
\finitem

\bibitem[BK]{bk}
\bgroup\bf V.~Bangert\egroup{} and \bgroup\bf M.~Katz\egroup{}.
Riemannian manifolds with harmonic one-forms of constant norm. Preprint,
\vskip -0.8em
\noindent {\bf http://www.math.biu.ac.il/\-\~{}katzmik/\-publications.html}.
\finitem

\bibitem[Bes87]{besse}
\bgroup\bf A.~L. Besse\egroup{}.
\newblock {\em Einstein Manifolds}.
\newblock Springer-Verlag, 1987.
\finitem

\bibitem[ES64]{eesa}
\bgroup\bf J.~Eells\egroup{} and \bgroup\bf J.P. Sampson\egroup{}.
\newblock Harmonic mappings of riemannian manifolds.
\newblock {\em Am. J. Math.}, 86:109--160, 1964.
\finitem

\bibitem[Kot01]{kot}
\bgroup\bf D.~Kotschick\egroup{}.
\newblock On products of harmonic forms.
\newblock {\em Duke Math. J.}, 107(3):521--531, 2001.
\finitem


\bibitem[Leb02]{lebrun}
\bgroup\bf C.~LeBrun\egroup{}.
\newblock Hyperbolic manifolds, harmonic forms, and Seiberg-Witten invariants.
\newblock {\em Geometriae Dedicata}, 91:137--154, 2002.
\finitem

\bibitem[Lic69]{lich}
\bgroup\bf A.~Lichnerowicz\egroup{}.
\newblock Application harmoniques dans un tore.
\newblock {\em C.R. Acad. Sci. S{\'e}r. A}, 269:912--916, 1969.
\finitem

\bibitem[Nag01]{nagy}
\bgroup\bf P.A.~Nagy\egroup{}.
\newblock {\em Un principe de s\'eparation des variables pour le spectre du laplacien des formes diff\'erentielles et applications}.
\newblock Th\`ese de doctorat, universit\'e de Savoie, 2001.
\newblock {\bf http://www.unine.ch/\-math/\-personnel/\-equipes/nagy/PN.htm}
\finitem

\bibitem[OLS78]{olszak}
\bgroup\bf Z.~Olszak\egroup{}.
\newblock A note on almost K\"ahler manifolds.
\newblock {\em Bull. Acad. Polon. Sci}, XXVI:199--206, 1978.
\finitem

\bibitem[Oni93]{ency}
\bgroup\bf A.L. Onishchik\egroup{}, editor.
\newblock {\em Lie groups and Lie Algebras I}, volume~20 of {\em Encyclopaedia
  of Mathematical Sciences}.
\newblock Springer-Verlag, 1993.
\finitem

\bibitem[PS61]{palais}
\bgroup\bf R.~S. Palais\egroup{} and \bgroup\bf T.~E. Stewart\egroup{}.
\newblock Torus bundles over a torus.
\newblock {\em Proc. Amer. Math. Soc.}, 12:26--29, 1961.
\finitem

\bibitem[PT01]{paktak}
\bgroup\bf H.K. Pak\egroup{} and \bgroup\bf T.~Takahashi\egroup{}.
\newblock Harmonic forms in a compact contact manifold.
\newblock {\em Tohoku Math. Publ.}, 20:125--138, 2001.
\finitem

\bibitem[Spi79]{spivak}
\bgroup\bf M.~Spivak\egroup{}.
\newblock {\em A comprehensive Introduction to Differential Geometry}.
\newblock Publish or Perish, inc., 1979.
\finitem

\bibitem[Ver01]{vernicos}
\bgroup\bf C.~Vernicos\egroup{}.
\newblock {\em The Macroscopic Spectrum of Nilmanifolds with an Emphasis on the
Heisenberg Groups}.
\newblock {\bf arXiv:math.DG/0210393}
\finitem

\bibitem[War83]{warner}
\bgroup\bf F.W. Warner\egroup{}.
\newblock {\em Foundations of Differentiable Manifolds and Lie Groups}.
\newblock Number~94 in Graduate texts in Mathematics. Springer-Verlag, 1983.
\finitem
\bibliend

\bigskip

\hbox{
\vbox{\hsize 6cm  \parindent=0cm \parskip=0cm
{\bf Paul-Andi Nagy}\par
\sf Universit\'e de Neuch\^atel\par
Institut de Math\'ematiques\par
11, rue Emile Argand\par
CH-2007 Neuch\^atel\par
Switzerland\par
mail: \tt Paul.Nagy@unine.ch}

\vbox{\hsize 6.5cm  \parindent=0cm \parskip=0cm
{\bf Constantin Vernicos}\par
\sf Universit\'e de Neuch\^atel\par
Institut de Math\'ematiques\par
11, rue Emile Argand\par
CH-2007 Neuch\^atel\par
Switzerland\par
mail: \tt Constantin.Vernicos@unine.ch}
}

\end